\documentclass{article}
\usepackage{amsmath}
\numberwithin{equation}{section}
\usepackage{indentfirst}

\begin{document}

\newtheorem{thm}{Theorem}[section]
\newtheorem{cor}[thm]{Corollary}
\newtheorem{prop}[thm]{Proposition}
\newtheorem{conj}[thm]{Conjecture}
\newtheorem{lem}[thm]{Lemma}
\newtheorem{Def}[thm]{Definition}
\newtheorem{rem}[thm]{Remark}
\newtheorem{prob}[thm]{Problem}
\newtheorem{ex}{Example}[section]

\newcommand{\be}{\begin{equation}}
\newcommand{\ee}{\end{equation}}
\newcommand{\ben}{\begin{enumerate}}
\newcommand{\een}{\end{enumerate}}
\newcommand{\beq}{\begin{eqnarray}}
\newcommand{\eeq}{\end{eqnarray}}
\newcommand{\beqn}{\begin{eqnarray*}}
\newcommand{\eeqn}{\end{eqnarray*}}
\newcommand{\bei}{\begin{itemize}}
\newcommand{\eei}{\end{itemize}}

\newcommand{\pa}{{\partial}}
\newcommand{\V}{{\rm V}}
\newcommand{\R}{{\rm R}}
\newcommand{\e}{{\epsilon}}
\newcommand{\tomega}{\tilde{\omega}}
\newcommand{\tOmega}{\tilde{Omega}}
\newcommand{\tR}{\tilde{R}}
\newcommand{\tB}{\tilde{B}}
\newcommand{\tGamma}{\tilde{\Gamma}}
\newcommand{\fa}{f_{\alpha}}
\newcommand{\fb}{f_{\beta}}
\newcommand{\faa}{f_{\alpha\alpha}}
\newcommand{\faaa}{f_{\alpha\alpha\alpha}}
\newcommand{\fab}{f_{\alpha\beta}}
\newcommand{\fabb}{f_{\alpha\beta\beta}}
\newcommand{\fbb}{f_{\beta\beta}}
\newcommand{\fbbb}{f_{\beta\beta\beta}}
\newcommand{\faab}{f_{\alpha\alpha\beta}}

\newcommand{\pxi}{ {\pa \over \pa x^i}}
\newcommand{\pxj}{ {\pa \over \pa x^j}}
\newcommand{\pxk}{ {\pa \over \pa x^k}}
\newcommand{\pyi}{ {\pa \over \pa y^i}}
\newcommand{\pyj}{ {\pa \over \pa y^j}}
\newcommand{\pyk}{ {\pa \over \pa y^k}}
\newcommand{\dxi}{{\delta \over \delta x^i}}
\newcommand{\dxj}{{\delta \over \delta x^j}}
\newcommand{\dxk}{{\delta \over \delta x^k}}

\newcommand{\px}{{\pa \over \pa x}}
\newcommand{\py}{{\pa \over \pa y}}
\newcommand{\pt}{{\pa \over \pa t}}
\newcommand{\ps}{{\pa \over \pa s}}
\newcommand{\pvi}{{\pa \over \pa v^i}}
\newcommand{\ty}{\tilde{y}}
\newcommand{\bGamma}{\bar{\Gamma}}

\font\BBb=msbm10 at 12pt
\newcommand{\Bbb}[1]{\mbox{\BBb #1}}

\newcommand{\qed}{\hspace*{\fill}Q.E.D.}  

\title{A class of Randers metrics of scalar flag curvature}
\author{ Xinyue Cheng\footnote{supported by the National Natural Science Foundation of China (11871126) and the Science Foundation of Chongqing Normal University (17XLB022)}, Li Yin, Tingting Li}

\date{}

\maketitle

\begin{abstract}
One of the most important problems in Finsler geometry is to classify Finsler metrics of scalar flag curvature. In this paper, we study the classification problem of Randers metrics of scalar flag curvature. Under the condition that $\beta$ is a Killing 1-form, we obtain some important necessary conditions for  Randers metrics to be of scalar flag curvature. 
{\bf  MR(2000) Subject Classification:}  53B40,  53C60
\end{abstract}

\section{Introduction}

In Finsler geometry, the flag curvature is the natural analogue of sectional curvature in Riemannian geometry and is an important Riemannian geometric quantity. The flag curvature characterizes the shape of Finsler spaces, which was firstly introduced by L. Berwald (\cite{B2,B1}). For a Finsler manifold $(M, F)$, the flag curvature ${\bf K=K} (x,y,P)$ of $F$ is a function of ``flag" $P\in T_{x}M $ and ``flagpole" $y\in T_{x}M $ at $x$ with $y\in P $. A Finsler metric $F$ is said to be of {\it scalar flag curvature} if for any non-zero vector $y \in T_xM, {\bf K=K}(x,y)$ is independent of $P$ containing $y \in T_xM$ (that is, the flag curvature is just a scalar function on the slit tangent bundle). In particular, a Finsler metric $F$ is said to be of {\it weakly isotropic flag curvature} if its flag curvature is a scalar function on $TM$ in the  form  ${\bf K}=\frac{3c_{x^m}y^m}{F(x,y)}+\sigma(x)$, where $c=c(x)$ and $\sigma =\sigma (x)$ are scalar functions on $M$.  A Finsler metric $F$ is of {\it constant flag curvature} if ${\bf K=}constant$. In 1929, Berwald proved that a projectively flat Finsler metric must be of scalar flag curvature (\cite{Ber3}).  However, we can find infinitely many Finsler metrics of scalar flag curvature, which are not projective flat. Therefore, one of the most important problems in Finsler geometry is to study and characterize Finsler metrics of scalar flag curvature. As we known, every two-dimensional Finsler metric is of scalar flag curvature. It is then a natural problem to understand Finsler metrics of scalar flag curvature in dimension $n \geq 3$.

Randers metrics are the simplest non-Riemannian Finsler metrics that can be expressed in the following special form
\[
F=\alpha+\beta,
\]
where $\alpha=\sqrt{a_{ij}y^iy^j}$ is a Riemannian metric and $\beta=b_iy^i$ is a non-zero 1-form on $M$ such that the norm of $\beta$ with respect to $\alpha$ satisfies that $||\beta||_{\alpha}<1$. Randers metrics were introduced by physicist G. Randers in 1941 from the standpoint of general relativity (\cite{R1}). Later on, these metrics were used in the theory of the electron microscope by R. S. Ingarden in 1957, who first named them Randers metrics. Randers metrics form an important class of Finsler metrics with a strong presence in both the theory and applications of Finsler geometry, and studying Randers metrics is an important step to understand general Finsler metrics. It is an important and fundamental problem in Finsler geometry to characterize and classify Randers metrics of scalar flag curvature.
Matsumoto and Shimada characterized Randers metrics of constant flag curvature (\cite{MaSh}). In 2003, Z. Shen classified locally projectively flat Randers metrics of constant flag curvature (\cite{Z3}). Further, D. Bao, Z. Shen and C. Robles showed that a Randers metric $F$ with navigation data $(h,W)$  is  of constant flag curvature if and only if $h$ is a Riemannian metric of constant sectional curvature and $W$ is homothetic with respect to $h$. Based on this, they classified Randers metrics of constant flag curvature via navigation data (\cite{BSR}). Furthermore, the first author, X. Mo and Z. Shen  classified projectively flat Randers metrics with isotropic S-curvature (\cite{CMS}). In this case, the metrics must be of weakly isotropic flag curvature  ${\bf K}=\frac{3c_{x^m}y^m}{F(x,y)}+\sigma(x)$.  Later,  the first author and Z. Shen classified  completely Randers metrics of weakly isotropic flag curvature on the manifolds of dimension $n\geq 3$ (\cite{CS1}).

For the current research on the classification problem of Randers metrics of scalar flag curvature, it is very important to find and construct some new examples. B. Chen and L. Zhao  constructed a new class of Randers metrics of scalar flag curvature (\cite{CZ}). Further, Z. Shen and Q. Xia proved firstly the following theorem:
\begin{quote}
{\bf Theorem A}(\cite{SX}) Let $\bar{F}=\alpha +\bar{\beta}$ be a Randers metric of weakly isotropic flag curvature on a manifold $M$ with navigation data $(\bar{h}, \bar{W})$ and $V$ be a homothetic vector field on $(M,\bar{F})$ with dilation $c$. Let $F=\alpha +\beta$ be a Randers metric on $M$ such that $\eta :=\beta -\bar{\beta}$ is closed. Assume that $\eta$ satisfies the following
\[
V^{i}\eta _{j;i}+\eta ^{i}V_{i;j}= 2c\eta _{j},
\]
where ``;" denotes the covariant derivative with respect to the Levi-Civita connection of $\bar{h}$ and $V_{i}:=\bar{h}_{ij}V^{j}, \ \eta ^{i}:=\bar{h}^{ij}\eta _{j}$. Then
\bei
\item $V$ is also homothetic vector field with respect to $F$.
\item With this $V$ such that $F(x, -V_{x})<1$, define $\tilde{F}=\tilde{\alpha}+\tilde{\beta}$ by
\be
F\left(x, \frac{y}{\tilde{F}}-V\right)=1. \label{barF}
\ee
Then $\tilde{F}$ is of scalar flag curvature and its flag curvature is related to that of $F$ by
\[
{\bf K}_{\tilde{F}}(x,y)={\bf K}_{F}(x,\tilde{y})-c^2 ,
\]
where $\tilde{y}:=y-\tilde{F}(x,y)V$.
\eei
\end{quote}
Then, based on Theorem A,  they found some new examples of Randers metrics of scalar flag curvature (\cite{SX}).

Soon after, Q. Xia studied Randers metrics of Douglas type and  proved the following theorem:
\begin{quote}
{\bf Theorem B} (\cite{Xia}) Let $F=\alpha +\beta$ be a Randers metric of Douglas type on an $n$-dimensional manifold $M (n\geq 3)$ and $V$ be a conformal vector field on $(M,F)$ with a conformal factor $c(x)$. Let $\tilde{F}=\tilde{\alpha}+\tilde{\beta}$ be a Randers metric defined from $(F,V)$ given by (\ref{barF}).
Then each two of the following imply the third one:
\ben
\item[{\rm (1)}] $F$ is of scalar flag curvature.
\item[{\rm (2)}] $\tilde{F}$ is of scalar flag curvature.
\item[{\rm (3)}] $V$ is homothetic or $\beta =0$.
\een
\end{quote}
Then, by using a special Randers metric $F=|y|+\langle x,y\rangle$ and a conformal vector field $V=xQ$ with $|xQ|<1$ on $R^n$, where $Q$ is an antisymmetric matrix, Q. Xia constructed a new Randers metric $\tilde{F}=\tilde{\alpha}+\tilde{\beta}$ by solving (\ref{barF}) such that $\tilde{F}$ is of scalar flag curvature but $\tilde{F}$ is not locally projectively flat. In this case, $V$ is actually a Killing vector field with respect to $F$ (\cite{Xia}).

However, although a significant progress has been made in studying Randers metrics of scalar flag curvature as mentioned above, it is still an important open problem in Finsler geometry to characterize and classify Randers metrics of scalar flag curvature.

For a 1-form $\beta=b_iy^i$ on $M$, we say that $\beta$ is a Killing 1-form with respect to a Riemannian metric $\alpha=\sqrt{a_{ij}y^iy^j}$ if it satisfies
\[
b_{i;j}+b_{j;i}=0,
\]
where $``\ ; "$ is the covariant derivative with respect to Levi-Civita connection of $\alpha$.
In this paper, we study classification problem of Randers metrics of scalar flag curvature. Under the condition that $\beta$ is a Killing 1-form, we obtain some necessary conditions for Randers metrics to be of scaler flag curvature. Firstly, we have the following theorem.

\begin{thm}\label{Randers metric1} \ Let $F=\alpha+\beta$ be a Randers metric on an $n$-dimensional manifold $M$. Suppose that $\beta$ is a Killing 1-form with respect to $\alpha$.  If $F$ is of scalar flag curvature ${\bf K=K}(x,y)$ , then $\alpha$ and $\beta$ satisfy the following equations
\beq
&& s^m_{\ 0;m}=-(n-1)b^{-2}c(x)\beta, \label{main11}\\
&& t_{00}+s_{0;0}=c(x)(\alpha^2-b^{-2}\beta^2),\label{main12}
\eeq
where $t_{ij}=s_{im}s^m_{\ j}$, $t_{00}=t_{ij}y^iy^j$ and $b:=||\beta_{x}||_{\alpha}$ denotes the norm of $\beta$ with respect to $\alpha$ and  $c=c(x)$ is a scalar function on $M$.
\end{thm}

Further, based on Theorem \ref{Randers metric1} and Proposition \ref{Randers metric6} below, we can obtain the following theorem.

\begin{thm}\label{Randers metric2} \ Let $F=\alpha+\beta$ be a Randers metric on an $n$-dimensional manifold $M$. Suppose that $\beta$ is a Killing 1-form with respect to $\alpha$. If $F$ is of scalar flag curvature ${\bf K=K}(x,y)$, then $\alpha$ and $\beta$ satisfy the following equation
\be
t_{0}=-\frac{n-1}{n+1}(\lambda+cb^{-2})\beta,\label{thm2}
\ee
where $t_i=b^jt_{ji}$, $t_0=t_iy^i$, $c=c(x)$ and $\lambda=\lambda(x)$ are scalar functions on  $M$.
\end{thm}

As an application of Theorem \ref{Randers metric1} and Theorem \ref{Randers metric2}, we study Randers metric on $S^3$ constructed by D. Bao and Z. Shen in \cite{Bao-Shen}. Bao-Shen have proved that this Randers metric is of constant flag curvature $K>1$. We will show that this metric also satisfies (\ref{main11}), (\ref{main12}) and (\ref{thm2}) (see Section \ref{example}).  Besides, Theorem \ref{Randers metric1} and Theorem \ref{Randers metric2} tell us, if a Randers metric $F=\alpha +\beta$ with Killing 1-form $\beta$  does not satisfy one of (\ref{main11}), (\ref{main12}) and (\ref{thm2}), then $F$ is not of scalar flag curvature.

\section{Preliminaries}

Let $F$ be a Finsler metric on an $n$-dimensional manifold $M$. The geodesics of  $F=F(x,y)$ are characterized by
\[
\frac{d^2 x^i}{d t^2}+2 G^i(x(t),\frac{dx}{dt}(t))=0,
\]
where $G^i$ are  the geodesic coefficients of $F$.

For any $x\in M$ and $y\in T_{x}M\backslash\{0\}$, the {\it Riemann curvature} ${\bf R}_{y} = R^i_{\ k} \frac{\pa }{\pa x^i} \otimes dx^k$  is defined by
\[
R^i_{\ k} = 2 \frac{\pa G^i}{\pa x^k} -\frac{\pa^2 G^i}{\pa x^m\pa y^k}y^m + 2 G^m \frac{\pa^2 G^i}{\pa y^m \pa y^k} -\frac{\pa G^i}{\pa y^m} \frac{\pa G^m}{\pa y^k}. \label{Rik}
\]
The trace of the Riemann curvature is {\it Ricci curvature} {\bf Ric}, which is defined by
\be
{\bf Ric}=R^m_{\ m}.
\ee
Obviously, the Ricci curvature is a positive homogenous function of degree two in $y$. The {\it Ricci tensor} is defined by
\[
{\bf Ric}_{ij}:=\frac12{\bf Ric}_{y^i y^j}.
\]

The flag curvature of Finsler manifold $(M,F)$ is the function ${\bf K}={\bf K}(x,y, P)$ of a two-dimensional plane called ``flag" $P \subset T_{x}M$ and a ``flagpole" $y\in P \backslash\{0\}$ defined by
$$
{\rm \bf K}(x,y, P):=\frac{g_{y}({\rm \bf R}_{y}(u),u)}{g_{y}(y,y)g_{y}(u,u)-\big[g_{y}(u,y)\big]^{2}},
$$
where $P =Span\{y, u\}$.

It is known that $F$ is of scalar flag curvature if and only if, in a standard local coordinate system£¬ we have (\cite{chernshen})
\be
R^i_{\ k} ={\bf{K}}(x,y) \left({F^2 \delta^{i}_{\ k}-F F_{y^{k}} y^i}\right),
\ee
which means
\be
{\bf{Ric}}=(n-1){\bf{K}}(x,y)F^2.    \label{2.1}
\ee

In Finsler geometry, there are some important quantities which all vanish for Riemannian metrics. Hence they are said to be {\it non-Riemannian}.  Let $\{{\bf b}_{i}\}$ be a basis for $T_{x}M$ and $\{\omega ^{i}\}$ be the basis for $T^{*}_{x}M$ dual to $\{{\bf b}_{i}\}$. Define the
Busemann-Hausdorff volume form by
\[
d V _{BH}:=\sigma _{BH}(x)\omega ^{1}\wedge \cdots \wedge\omega ^{n},
\]
where
\[
\sigma _{BH}(x):=\frac{Vol\big({\bf B}^{n}(1)\big)}{Vol\left\{(y^{i})\in R^{n}|F(x,y^{i}{\bf b}_{i})<1 \right\}}.
\]
Here $Vol\{ \cdot \}$ denotes the Euclidean volume function on subsets in ${\bf R}^{n}$ and ${\bf B}^{n}(1)$ denotes the unit ball in ${\bf R}^{n}$.

Define
$$
\tau (x,y):=\ln \left[\frac{\sqrt{det(g_{ij}(x,y))}}{ \sigma _{BH}(x)}\right].
$$
$\tau =\tau (x,y)$ is well-defined, which is called the {\it distortion} of $F$. The distortion $\tau$ characterizes the geometry of tangent space $(T_{x}M, F_{x})$.
It is natural to study the rate of change of the distortion along geodesics. For a vector $y\in T_{x}M\setminus\{0\}$, let $\sigma =\sigma (x)$ be the geodesic with $\sigma (0)=x$ and $\dot{\sigma}(0)=y$.
Put
\[
{\bf S}(x,y):=\frac{d}{dt}\left[\tau\left(\sigma (t),\dot{\sigma}(t) \right) \right]|_{t=0}.
\]
Equivalently,
$$
{\bf S}(x,y):=\tau _{|m}(x,y)y^{m},
$$
where $``|"$ denotes  the horizontal covariant derivative with respect to Chern connection of $F$. ${\bf S}={\bf S}(x,y)$ is called the {\it S-curvature} of Finsler metric $F$.
The S-curvature ${\bf S}(x,y)$ measures the rate of  change of $(T_{x}M, F_{x})$ in the direction $y\in T_{x}M$ (\cite{chernshen,Shen2}). In particular, ${\bf S}=0$ for Riemannian metrics.

There are other several important non-Riemannian quantities in Finsler geometry. Let
\[
I_{i}(x,y):=\frac{\pa \tau}{\pa y^i}(x,y)=\frac{1}{2}g^{jk}(x,y)\frac{\pa g_{jk}}{\pa y^i}(x,y).
\]
The tensor ${\bf I}_y:=I_i(x,y)dx^i$ is called the {\it mean Cartan torsion} of $F$. According to Deicke's theorem, a Finsler metric $F$ is Riemannian if and only if ${\bf I}_y=0$.

Another important non-Riemannian quantity is the {\it Landsberg curvature} ${\bf L}=L_{ijk}(x,y)dx^i\otimes dx^y\otimes dx^k$  defined by
\[
L_{ijk}(x,y):=-\frac{1}{2}F F_{y^{m}} \frac{\pa^3 G^m}{\pa y^i \pa y^j \pa y^k}.
\]
A Finsler metric $F$ is called {\it the Landsberg metric} if the Landsberg curvature vanishes, i.e. ${\bf L}=0$. Further, let $J_k:=g^{ij}L_{ijk}$. Then ${\bf J}=J_k dx^k$  is called the {\it mean Landsberg curvature}.
A Finsler metric $F$ is said to be {\it weakly Landsbergian} if ${\bf J}=0$.

Now we assume that  Finsler metric $F$ is of scalar flag curvature ${\bf K}={\bf K}(x,y)$. Firstly, we have the following identity (\cite{CMS,chernshen})
\be
{\bf S}_{\cdot i|m}y^m-{\bf S}_{|i}=-\frac{n+1}{3}{\bf K}_{\cdot i}F^2 , \label{2.2}
\ee
where ${{\bf K}_{\cdot i}:={\bf K}_{y^i}}$.
Further, we have following important identity which connects the mean Cartan torsion $I_i$ and the mean Landsberg curvature $J_i$ with the flag curvature  (\cite{chernshen})
\be
J_{i|m}y^m+{\bf K}F^2I_i=-\frac{n+1}{3}F^2{\bf K}_{\cdot i}. \label{2.3}
\ee

In the following, let us consider  $(\alpha, \beta)$-metrics in the form  $F=\alpha \phi (\beta/\alpha)$, where $\alpha=\sqrt{a_{ij} y^i y^j}$ is a Riemannian metric and  $\beta=b_i y^i$ is a  1-form on a manifold.
Let

\[
r_{ij}:=\frac{1}{2}(b_{i;j}+b_{j;i}),\ s_{ij}:=\frac{1}{2}(b_{i;j}-b_{j;i}),
\]
\[
r^i_{\ j}:=a^{ik}r_{kj},\ s^i_{\ j}:=a^{ik}s_{kj}, \  r_j:=b^i r_{ij},\ s_j:=b^i s_{ij},
\]
\[
q_{ij}:=r_{im} s^{m}_{\ j},\ t_{ij}:=s_{im} s^{m}_{\ j}, \ q_j:=b^i q_{ij},\  t_j:=b^i t_{ij},
\]
where $b^i:=a^{ij}b_j$, $(a^{ij})$ is the inverse of $(a_{ij})$ and `` ; " denotes the covariant derivative with respect to Levi-Civita connection of  $\alpha$.  We use  $(a^{ij})$ and $(a_{ij})$ to raise or lower the indices of a tensor.

For an $(\alpha, \beta)$-metric $F=\alpha \phi(s),\ s=\beta / \alpha$, we have the formulas for the mean Landsberg curvature (\cite{LiSh1}) and the mean Cartan torsion (\cite{CWW}) of $F$ as follows

\begin{align}
I_i=&-\frac{\Phi (\phi-s\phi ')h_i}{2 \Delta \phi \alpha^2},\label{2.4}\\
J_i=&- \frac{1}{2 \Delta \alpha^4}\bigg\{\frac{2 \alpha^2}{b^2-s^2}\Big[\frac{\Phi}{\Delta}+(n+1)(Q+s Q')\Big](s_0+r_0)h_i\nonumber\\
&+\frac{\alpha}{b^2-s^2}\Big[\Psi_1+s \frac{\Phi}{\Delta}\Big](r_{00}-2 \alpha Q s_0)h_i+\alpha\Big[-\alpha Q's_0h_i+\alpha Q (\alpha^2 s_i-y_is_0)\nonumber\\
&+\alpha^2 \Delta s_{i0}+\alpha^2(r_{i0}-2 \alpha Q s_i)-(r_{00}-2\alpha Q s_0)y_i\Big]\frac{\Phi}{\Delta} \bigg\}.\label{2.5}
\end{align}
Besides, we also have (\cite{LiSh1})
\be
\bar{J}:=J_i b^i=-\frac{1}{2 \Delta \alpha^2}\left\{\Psi_1(r_{00}-2 \alpha Q s_0)+\alpha\Psi_2(r_0+s_0) \right\},\label{2.6}
\ee
where
\beqn
&& h_i:=\alpha b_i-s y_i,\ \ Q:=\frac{\phi'}{\phi-s \phi'},\ \ \Delta:=1+s Q+(b^2-s^2)Q', \\
&&\Phi:=-(Q-s Q')\left\{ n\Delta+1+s Q\right\}-(b^2-s^2)(1+s Q)Q'' , \\
&&\Psi_1:=\sqrt{b^2-s^2} \Delta^{\frac{1}{2}}\left[\frac{\sqrt{b^2-s^2}\Phi}{\Delta^{\frac{3}{2}}}\right]', \ \Psi_2:=2 (n+1)(Q-s Q')+3 \frac{\Phi}{\Delta}.
\eeqn

\section{Some properties of Randers metrics of scalar flag curvature}\label{section3}

Let $F=\alpha +\beta$ be a Randers metric on an $n$-dimensional manifold $M$ and $\beta$ be a Killing 1-form (i.e., $r_{ij}=0$) with respect to $\alpha$. In this case,  $\phi=1+s$ and we have
\beqn
&&  Q=1,\ \  \Delta=1+s,\ \  \Phi=-(n+1)(1+s), \\
&&  \Psi_1=\frac{(n+1)(2s+s^2+b^2)}{2 (1+s)}, \ \  \Psi_2=-(n+1),    \\
&&\Psi_1'=\frac{(n+1)(2+2s-b^2+s^2)}{2 (1+s)^2},\\
&&\Psi_1+s\frac{\Phi}{\Delta}=\Psi_1-s(n+1)=\frac{(n+1)(b^2-s^2)}{2(1+s)},\\
&&2\Psi_1-\Psi_2=2\Psi_1+(n+1)=\frac{(n+1)(1+3s+b^2+s^2)}{1+s}.
\eeqn
From these and (\ref{2.4})-(\ref{2.6}), we get
\beq
J_i= \frac{n+1}{2 \alpha^2 (1+s)} \Big\{\frac{ s_0h_i}{1+s}+\alpha (1+s) s_{i0}-\alpha^2 s_i+ s_0y_i\Big\},\label{3.1}
\eeq
and
\beq
I_i=\frac{n+1}{2 \alpha^2 (1+s)}h_i,\ \ \ \
\bar{J}=\frac{n+1}{2\alpha(1+s)^2}(1+3s+b^2+s^2)s_0.\label{3.2}
\eeq

In the following, we always assume that Randers metric $F=\alpha+\beta$ is of scalar flag curvature. In this case, by the formulas (4.15) and (5.30) for Ricci curvature of $F$ in  \cite{chsh3}, we have
\be
{\bf Ric}=2 \alpha s^{m}_{\ 0;m}+(n-1)\left(\kappa \alpha^2 +t_{00}-\frac{2 \alpha^2}{F} t_{0}+\frac{3 \alpha^2 s_{0}^2}{F^2}+\frac{\alpha s_{0;0}}{F}\right), \label{3.3}
\ee
where $\kappa=\kappa(x)$ is a scalar function on $M$.  Obviously, by (\ref{2.1}) and (\ref{3.3}), we can get

\begin{lem}\label{Randers metric3} \ Let $F=\alpha+\beta$ be a Randers metric on an $n$-dimensional manifold $M$. If $F$ is of scalar flag curvature ${\bf K=K}(x,y)$ and  $\beta$ is a Killing 1-form with respect to $\alpha$, then
\be
2 \alpha s^{m}_{\ 0;m}+(n-1)(\kappa \alpha^2 +t_{00}+\Xi)=(n-1){\bf K}F^2,\label{lem3.1}
\ee
where
\beq
\Xi:=-\frac{2 \alpha^2 t_{0}}{F}+\frac{3 \alpha^2 s_{0}^2}{F^2}+\frac{\alpha s_{0;0}}{F}.
\eeq

\end{lem}

Further, contracting the (\ref{2.2}) by $b^i$ yields the following equation
\beq
{\bf S}_{\cdot i|m}y^m b^i-{\bf S}_{|i} b^i=-\frac{n+1}{3}F^2{\bf K}_{\cdot i}b^{i}.\label{3.5}
\eeq
Let $G^{i}$ and $\bar{G^{i}}$ denote the geodesic coefficients of $F$ and $\alpha$ respectively. The horizontal covariant derivatives ${\bf S}_{\cdot i|m}$ and ${\bf S}_{\cdot i;m}$ of ${\bf S}_{\cdot i}$ with respect to $F$ and $\alpha$ respectively are given by
\beq
&&{\bf S}_{\cdot i|m}=\frac{\pa {\bf S}_{\cdot i} }{\pa x^m}-{\bf S}_{\cdot l}\Gamma^l_{\ im}-\frac{\pa {\bf S}_{\cdot i} }{\pa y^l}N^l_m, \label{Sim1} \\
&&{\bf S}_{\cdot i;m}=\frac{\pa {\bf S}_{\cdot i} }{\pa x^m}-{\bf S}_{\cdot l}{\bar\Gamma}^l_{\ im}-\frac{\pa {\bf S}_{\cdot i} }{\pa y^l}{\bar N}^l_m,\label{Sim2}
\eeq
where $\Gamma^l_{\ im}=\frac{\pa^2 G^l }{\pa y^i y^m}-L^{l}_{\ im}$, \ $N^l_m:=\frac{\pa^2 G^l }{\pa y^m}, \ L^{l}_{\ im}:= g^{jl}L_{jim}$ and ${\bar\Gamma}^l_{\ im}:=\frac{\pa^2 {\bar G}^l }{\pa y^i y^m} $, ${\bar N}^l_m:=\frac{\pa^2 {\bar G}^l }{\pa y^m}$.
From (\ref{Sim1}) and (\ref{Sim2}), we have
\beq
{\bf S}_{\cdot i|m}y^m &=& \Big\{{\bf S}_{\cdot i;m}-{\bf S}_{\cdot l}(\Gamma^l_{\ im}-{\bar\Gamma}^l_{\ im})-\frac{\pa {\bf S}_{.i} }{\pa y^l}(N^l_m-{\bar N}^l_m)\Big\}y^m\nonumber\\
&=& {\bf S}_{\cdot i;m}y^m-{\bf S}_{\cdot l}(N^l_i-{\bar N}^l_i)-2\frac{\pa {\bf S}_{\cdot i} }{\pa y^l}(G^l-{\bar G}^l).\label{3.6}
\eeq
Similarly, we have
\be
{\bf S}_{|i}={\bf S}_{;i}-{\bf S}_{\cdot l}(N^l_i-{\bar N}^l_i).\label{3.7}
\ee
Substituting (\ref{3.6}) (\ref{3.7}) into (\ref{3.5}) yields
\be
{\bf S}_{.i;m}y^m b^i-2\frac{\pa {\bf S}_{.i} }{\pa y^l}(G^l-{\bar G}^l)b^i-{\bf S}_{;i}b^i=-\frac{n+1}{3}F^2{\bf K}_{.i}b^i.\label{3.8}
\ee

 As we known,  the relationship between the geodesic coefficients $G^i$  and ${\bar G}^i$ is as follows (\cite{chsh3,chernshen})
\be
G^i={\bar G}^i+\alpha s^i_{\ 0}-\frac{\alpha s_0 y^i}{F}.\label{3.9}
\ee
From this, one obtains
\beq
N^i_{\ j}& =& {\bar N}^i_{\ j}+\alpha^{-1}y_j s^i_{\ 0}+\alpha s^i_{\ j}-\frac{ s_0 y^i(y_j+\alpha b_j)}{F^2}\nonumber \\
           && -\frac{\alpha^{-1}y_j s_{0}y^i+\alpha y^is_j+\alpha s_0 \delta^i_{\ j}}{F}.\label{3.9.1}
\eeq
Besides, the formula for S-curvature of $F$ is given by the following (\cite{chernshen})
\be
{\bf S}=-(n+1)\left\{{\frac{\alpha s_0}{F}+\rho_0}\right\}.\label{3.4}
\ee
where $\rho:=ln\sqrt{1-||\beta||^2_{\alpha}}$, $\rho_i=\rho_{x^i}$ and $\rho_0=\rho_i y^i$.

By a series of direct computations, we have the following
\beq
&&  F_{\cdot j}=\alpha^{-1}y_j+b_j,\ \ F_{.j.k}=\alpha^{-1}a_{jk}-\alpha^{-3}y_j y_k, \ \ F_{.i}b^i=s+b^2,  \label{fijk}  \\
&&  F_{;i}=s_{0i}, \ \ F_{;i}y^i=0, \ \  F_{;i}b^i=-s_0, \\
&&  F_{\cdot j}s^j_{\ 0}=s_0, \ \  F_{\cdot k\cdot j}s^j_{\ 0}=\alpha^{-1}s_{k0}. \label{3.10}
\eeq
From (\ref{3.9}) and (\ref{3.4})-(\ref{3.10}),  one obtains
\begin{align}
&{\bf S}_{;i}b^i=(n+1)\Big\{-\frac{b^i s_{0;i}}{1+s}-\frac{ s_0^2}{\alpha (1+s)^2}-b^i\rho_{0;i}\Big\},\label{3.11}\\
&{\bf S}_{.i;m}y^m b^i=(n+1)\Big\{\frac{(b^2-s^2)s_{0;0}}{\alpha (1+s)^2}+\frac{ s_0^2}{\alpha (1+s)^2}-\frac{ b^i s_{i;0}}{1+s}-b^i\rho_{i;0}\Big\},\label{3.12}\\
&\frac{\pa {\bf S}_{.i} }{\pa y^l}(G^l-{\bar G}^l)b^i=(n+1)\Big\{-\frac{2 s_0^2(s+b^2)}{\alpha (1+s)^3}+\frac{ t_0(b^2-s^2)}{(1+s)^2}\Big\}.\label{3.13}
\end{align}
Plugging (\ref{3.11})-(\ref{3.13}) into (\ref{3.8}), we obtain the following lemma.

\begin{lem}\label{Randers metric4} \ Let $F=\alpha+\beta$ be a Randers metric on an $n$-dimensional manifold $M$. If $F$ is of scalar flag curvature ${\bf K=K}(x,y)$ and  $\beta$ is a Killing 1-form with respect to $\alpha$,  then
\beq
&&(n+1)\Big\{\alpha^{-1}\left[ \frac{s_{0;0}(b^2-s^2)}{(1+s)^2}+ \frac{2s_0^2(1+3s+2b^2)}{(1+s)^3}\right]+ b^i(\rho_{0;i}-\rho_{i;0})\nonumber\\
&& +\frac{2t_{0}(s^2-b^2)}{(1+s)^2}+\frac{b^i (s_{0;i}-s_{i;0})}{1+s}\Big\}=-\frac{n+1}{3}\alpha^2(1+s)^2{\bf K}_{.i}b^i .\label{lem3.2}
\eeq
\end{lem}

On the other hand, contracting (\ref{2.3}) by $b^i$ yields the following equation
\be
J_{i|m}y^m b^i+{\bf K}F^2I_ib^i=-\frac{n+1}{3}F^2{\bf K}_{.i}b^i. \label{3.14}
\ee
Similar to the methods to get Lemma \ref{Randers metric4}, we have
\beq
J_{i|m}y^m=J_{i;m}y^m-J_l \frac{\pa( G^l-{\bar G}^l)}{\pa y^i}-2\frac{\pa J_{i}}{\pa y^l}( G^l-{\bar G}^l).\label{3.15}
\eeq
Substituting (\ref{3.15}) into (\ref{3.14}) yields
\beq
&&{\bar J}_{;0}-J_i s^i_{\ 0}-J_l \frac{\pa( G^l-{\bar G}^l)}{\pa y^i}b^i-2\frac{\pa {\bar J}}{\pa y^l}( G^l-{\bar G}^l)+{\bf K} F^2 I_{i} b^i  \nonumber\\
&& =-\frac{n+1}{3}F^2 {\bf K}_{\cdot i} b^i,\label{3.16}
\eeq
where ${\bar J}=J_i b^i$.

It is easy to see that
\be
h_i s^i_{\ 0}=\alpha s_0,\ \ h_is^i=-s s_0, \ \ h_i y^i=0,\ \ h_i b^i=\alpha (b^2-s^2). \label{his}
\ee
From  (\ref{3.1}), (\ref{3.2}), (\ref{3.9}) and (\ref{3.9.1}), (\ref{his}),  by a series of direct computations, one obtains the following

\begin{align}
&I_ib^i=\frac{(n+1)(b^2-s^2)}{2\alpha (1+s)},\label{3.17}\\
&{\bar J}_{;0}=\frac{(n+1) \big[(1+3s+b^2+s^2)s_{0;0}+2 s_0^2\big]}{2 \alpha (1+s)^2},\\
&J_i s^i_{0}=\frac{n+1}{2 (1+s) \alpha}\left\{\frac{s_0^2}{1+s}-(1+s)t_{00}-\alpha t_0\right\},\\
&\frac{\pa {\bar J}}{\pa y^l}(G^l-{\bar G}^l)=(n+1)\left\{\frac{1+3s+b^2+s^2}{2 (1+s)^2}t_0+\frac{(1-s-2b^2)s_0^2}{2 \alpha (1+s)^3}\right\},\\
&J_l \frac{\pa (G^l-{\bar G}^l)}{\pa y^i}b^i=\frac{n+1}{2 \alpha (1+s)}\Big\{ -\alpha (1+2s) t_0-s (1+s)t_{00}  -\alpha^2 t_{i}b^i   \nonumber\\
&\ \ \ \ \ \ \ \ \ \ \ \ \ \ \ \ \ \ \ \ \ \ \ \ -\frac{(2+3s+b^2) s_0^2}{(1+s)^2}\Big\}.\label{3.18}
\end{align}
Substituting (\ref{3.17})-(\ref{3.18}) into (\ref{3.16}) , we have the following lemma.

\begin{lem}\label{Randers metric5} \ Let $F=\alpha+\beta$ be a Randers metric on an $n$-dimensional manifold $M$. If $F$ is of scalar flag curvature ${\bf K=K}(x,y)$ and   $\beta$ is a Killing 1-form with respect to $\alpha$,  then
\beq
&&\Big\{\alpha^{-1}\left[\frac{s_{0;0}(1+3s+b^2+s^2)}{2(1+s)^2}+\frac{t_{00}(1+s)}{2}+\frac{s_0^2(1+6s+5b^2)}{2(1+s)^3}\right] \nonumber \\
&& -\frac{t_0(s+b^2)}{(1+s)^2} +\frac{\alpha t_i b^i}{2(1+s)}+\frac{\alpha(1+s)(b^2-s^2){\bf K} }{2}\Big\}(n+1)\nonumber\\
&& =-\frac{n+1}{3}\alpha^2(1+s)^2{\bf K}_{.i}b^i.  \label{lem3.3}
\eeq
\end{lem}

Lemmas \ref{Randers metric3}, \ref{Randers metric4} and \ref{Randers metric5} are fundamental for the proofs of our main theorems.

\section{Proofs of The Theorems}

In this section, we will prove Theorem \ref{Randers metric1} and \ref{Randers metric2}. Firstly, we will prove the following proposition which plays an important role in our proofs.

\begin{prop}\label{Randers metric6}  \ Let $\alpha=\sqrt{a_{ij}y^iy^j}$ be a Riemann metric and $\beta=b_iy^i$ be a 1-form on a manifold. Then
\be
s^i_{\ 0;i}={^{\alpha}\bf{Ric}}_{0b}+r^i_{\ i;0}-r^i_{\ 0;i},\label{4.00}
\ee
where $`` ; "$ denotes the covariant derivative with respect to Levi-Civita connection of  $\alpha$, ${^{\alpha}\bf{Ric}}$ is the Ricci curvature of $\alpha$ and ${^{\alpha}\bf{Ric}}_{0b}:=(\frac12 {^{\alpha}\bf {Ric}})_{y^i y^j}y^ib^j$.
\end{prop}
{\it \bf Proof }: By the definition of $r_{ij}$, we have
\beq
b_{i;k;j}+b_{k;i;j}=2r_{ik;j},  \label{4.01}\\
-b_{k;j;i}-b_{j;k;i}= -2r_{kj;i}.
\eeq
On the other hand, by Ricci identity, we get
\beq
&& b_{i;j;k}-b_{i;k;j}= {^{\alpha}R}_{imjk}b^m,\\
&& b_{k;j;i}-b_{k;i;j}=-{^{\alpha}R}_{kmij}b^m,\\
&& b_{j;k;i}-b_{j;i;k}= {^{\alpha}R}_{jmki}b^m. \label{4.02}
\eeq
Here, ${^{\alpha}R}_{imjk}$ denote the Riemann curvature tensor of $\alpha$. Adding all the equations from (\ref{4.01}) to (\ref{4.02}), we can get
\be
s_{ij;k}=- ^{\alpha}R_{kmij}b^m+r_{ik;j}-r_{jk;i}.\label{4.03}
\ee
Contracting (\ref{4.03}) with $a^{pi}$ yields
\be
s^{p}_{\ j;k}= {^{\alpha}R}^{\ p} _{j\ kb}+r^p_{\ k;j}-r^{\ \ \ p}_{kj;},\label{4.04}
\ee
where ${^{\alpha}R}_{j \ km}^{\ p}:={^{\alpha}R}_{jikm}a^{ip}$ and ${^{\alpha}R}^{\ p} _{j\ kb}= {^{\alpha}R}^{\ p} _{j\ km}b^m$. Contracting (\ref{4.04}) with $y^j$ yields
\be
s^{p}_{\ 0;k}= {^{\alpha}R}^{\ p} _{0\ kb}+r^p_{\ k;0}-r^{\ \ \ p}_{k0;}.\label{4.05}
\ee
Letting $p=k$ in (\ref{4.05}), we obtain (\ref{4.00}).
\qed

\vskip 2mm

Now we are in the position to prove our main theorems.

\vskip 2mm

\noindent {\bf Proof of Theorem \ref{Randers metric1}}:  Let $F=\alpha+\beta$ be a Randers metric. Assume that $F$ is of scalar flag curvature ${\bf{K=K}}(x,y)$ and $\beta$ is a Killing 1-form. Then $(\ref{lem3.1})$, $(\ref{lem3.2})$ and $(\ref{lem3.3})$ all hold by Lemma \ref{Randers metric3}, Lemma \ref{Randers metric4} and Lemma \ref{Randers metric5}.

By $(\ref{lem3.3})-(\ref{lem3.2})+ (\ref{lem3.1}) \ \frac{b^2-s^2}{2 (n-1)\alpha (1+s)}$ and plugging $s=\frac{\beta}{\alpha}$ into the resulting equation, we can obtain
\[
\frac{1}{2\alpha (\alpha+\beta)}\left(\Xi_3 \alpha^3+\Xi_2 \alpha^2+\Xi_1 \alpha+\Xi_0\right)=0,
\]
which is equivalent to the following
\be
\Xi_3 \alpha^3+\Xi_2 \alpha^2+\Xi_1 \alpha+\Xi_0=0,\label{4.2}
\ee
where
\beq
&&\Xi_3=\kappa b^2+t_ib^i,\label{4.3}\\
&&\Xi_2=2b^i(\rho_{i;0}-\rho_{0;i}+s_{i;0}-s_{0;i})+\frac{2 b^2 s^m_{\ 0;m}}{n-1}, \label{Xi2} \\
&&\Xi_1=-\kappa \beta^2+[2 b^i(\rho_{i;0}-\rho_{0;i})-2t_0]\beta+t_{00}-3s_0^2+t_{00}b^2+s_{0;0},  \label{Xi1} \\
&&\Xi_0=2 \beta \left(s_{0;0}+t_{00}-\frac{\beta s^m_{\ 0;m}}{n-1}\right).  \label{4.4}
\eeq

Further, by the definition of $\rho$ and because of $r_{ij}=0$ , we have $\rho_{i}=\rho_{x^i}=-\frac{s_i}{1-b^2}$. Then, we can get
\[
\rho_{i;0}-\rho_{0;i}=\frac{(s_{0;i}-s_{i;0})}{1-b^2}.
\]

On the other hand, it is clear that $(b^2)_{i;j}=(b^2)_{j;i}$. By a direct computation, we get $s_{i;j}=s_{j;i}$, that is, $s_{i;0}=s_{0;i}$.
Hence, $(\ref{Xi2}), (\ref{Xi1})$ can be rewritten as
\beq
&&\Xi_2=\frac{2 b^2 s^m_{\ 0;m}}{n-1},\label{4.6}\\
&&\Xi_1=-\kappa \beta^2-2t_0\beta+t_{00}(1+b^2)-3s_0^2+s_{0;0}. \label{4.7}
\eeq
From (\ref{4.2}) we obtain the following fundamental equations
\beq
\Xi_2\alpha^2+\Xi_0=0,\label{4.9}\\
\Xi_3\alpha^2+\Xi_1=0.\label{4.10}
\eeq
Rewrite (\ref{4.9}) as
\be
\Xi_2\alpha^2=-\Xi_0.\label{4.11}
\ee
Since $\alpha^2$ is an irreducible polynomial in $y$, from (\ref{4.11}) and (\ref{4.4}), (\ref{4.6}), we have
\beq
s_{0;0}+t_{00}-\frac{ \beta s^m_{\ 0;m}}{n-1}=c(x)\alpha^2,\label{4.12}
\eeq
where $c=c(x)$ is a scalar function on $M$. Substituting (\ref{4.12}) into (\ref{4.11}) yields
\be
\frac{ b^2 s^m_{\ 0;m}}{n-1}=-c(x)\beta ,\label{4.13}
\ee
which means that (\ref{main11}) holds. Plugging (\ref{4.13}) into (\ref{4.12}), we get
\[
s_{0;0}+t_{00}=c(x)(\alpha^2-b^{-2}\beta^2),
\]
which is just (\ref{main12}). This completes the proof of Theorem \ref{Randers metric1}.
\qed

\vskip 2mm
Finally, let us give the proof of Theorem \ref{Randers metric2}.
\vskip 3mm
\noindent{\bf Proof of Theorem \ref{Randers metric2}}:  Let $F=\alpha+\beta$ be a Randers metric. Assume that $F$ is of scalar flag curvature ${\bf K=K}(x,y)$.  By (5.30) in \cite{chsh3}, we have the following
\[
^\alpha{\bf{Ric}}=(n-1)\lambda(x)\alpha^2 +(n+1)t_{00},
\]
where $\lambda:=\lambda(x)$ is a scalar function on $M$. Further, by a direct computation, we have
\be
{^\alpha{\bf{Ric}}}_{ij}=\left(\frac{1}{2} {{^\alpha{\bf{Ric}}}}\right)_{y^i y^j} =(n-1)\lambda(x)a_{ij} +(n+1)t_{ij}.\label{4.14}
\ee
From Proposition \ref{Randers metric6} and by (\ref{4.14}) and the assumption that $\beta$ is a Killing 1-form, we get
\be
s^m_{\ 0;m}=(n-1)\lambda(x)\beta+(n+1)t_{0}.\label{4.15}
\ee

On the other hand, by Theorem \ref{Randers metric1}, we have
\be
s^m_{\ 0;m}=-(n-1)cb^{-2}\beta.\label{4.16}
\ee
Comparing (\ref{4.15}) with (\ref{4.16}),  we obtain (\ref{thm2}).
\qed

\section{Example: Bao-Shen metric on $S^3$}\label{example}

In this section, we will prove that, for Bao-Shen metric on $S^3$ constructed in \cite{Bao-Shen}, $\alpha$ and $\beta$ satisfy (\ref{main11}), (\ref{main12}) and (\ref{thm2}).

We view $S^3$ as a Lie group and let $\eta ^1$, $\eta ^2$ and $\eta ^3$ be the standard right invariant 1-form on $S^3$ such that
\[
d\eta ^{1}=2\eta ^{2}\wedge \eta ^{3}, \ \ d\eta ^{2}=2\eta ^{3}\wedge \eta ^{1}, \ \ d\eta ^{3}=2\eta ^{1}\wedge \eta ^{2}.
\]
For a constant $K> 1$, write $\varepsilon :=\sqrt{K}, \ \delta :=\pm \sqrt{K-1}.$ Let $\omega ^{1}:=\varepsilon \eta ^{1}, \ \omega ^{2}:=\eta ^{2}, \ \omega ^{3}:=\eta ^{3}$. Then
\[
d\omega ^{1}=2\varepsilon \omega ^{2}\wedge \omega ^{3}, \ d\omega ^{2}=\frac{2}{\varepsilon}\omega ^{3}\wedge \omega ^{1}, \ d\omega ^{3}=\frac{2}{\varepsilon}\omega ^{1}\wedge \omega ^{2}.
\]
Consider the following Riemannian metric on $S^3$:
\be
\alpha :=\sqrt{\omega ^{1}\otimes \omega ^{1}+\omega ^{2}\otimes \omega ^{2}+\omega ^{3}\otimes \omega ^{3}}. \label{alpha}
\ee
The Levi-Civita connection forms $(\omega _{i}^{\ j})$ of $\alpha$ are given by $d\omega ^{i}=\omega ^{j}\wedge \omega _{j}^{\ i}$ and $\omega _{ij}= -\omega _{ji}$, where $\omega _{ij}:=\omega _{i}^{\ k}\delta _{kj}$. Concretely, we have
\[ \left(
\begin{array}{lcr}
\omega _{1}^{\ 1} & \omega _{1}^{\ 2} & \omega _{1}^{\ 3} \\
\omega _{2}^{\ 1} & \omega _{2}^{\ 2} & \omega _{2}^{\ 3} \\
\omega _{3}^{\ 1} & \omega _{3}^{\ 2} & \omega _{3}^{\ 3}
\end{array}
\right)=\left(
\begin{array}{ccc}
0 & -\varepsilon \omega ^{3} & \varepsilon \omega ^{2} \\
\varepsilon \omega ^{3} & 0 & \left(\varepsilon -\frac{2}{\varepsilon}\right)\omega ^{1} \\
-\varepsilon \omega ^{2} & \left(\frac{2}{\varepsilon}-\varepsilon \right) \omega ^{1} & 0
\end{array}
\right)
\]
Further, let
\be
\beta := \frac{\delta}{\varepsilon}\omega ^{1}. \label{beta}
\ee
Then we have
\[
b_{1}=\frac{\delta}{\varepsilon}, \ \ b_{2}=0, \ \ b_{3}=0.
\]
Obviously, this $\beta$ has constant length with respect to Riemannian metric $\alpha$:
\[
b:=\|\beta\|_{\alpha}=\frac{|\delta |}{\varepsilon}.
\]
By the Riemannian metric $\alpha$ defined by (\ref{alpha}) and the 1-form $\beta$ defined by (\ref{beta}), D. Bao and Z. Shen constructed the following family of Randers metrics on $S^3$:
\be
F(x,y)=\sqrt{(y^1)^2 +(y^2)^2 +(y^3)^2}\pm \sqrt{\frac{K-1}{K}}y^1 , \ \ \ \left(\frac{\delta}{\varepsilon}=\pm \sqrt{\frac{K-1}{K}}\right). \label{Bao-Shen}
\ee
Further, they have proved that this family of Randers metrics has constant flag curvature $K$ but nonzero Douglas tensor (see \cite{Bao-Shen}). This family of Randers metrics is the first example of Finsler metrics of constant flag curvature which are not projectively flat.

In the following, we carry out our calculations in the moving coframe $\{\omega ^{1}, \omega ^{2}, \omega ^{3}\}$ and the dual moving frame $\{e_{1}, e_{2}, e_{3}\}$. The covariant differentiation formulas of $\beta$ with respect to $\alpha$ are given by
\beqn
&&b_{i;j}=(db_{i}-b_{s}\omega _{i}^{\ s})(e_{j}), \\
&&b_{i;j;k}=(db_{i;j}-b_{s;j}\omega _{i}^{\ s}-b_{i;s}\omega _{j}^{\ s})(e_{k}).
\eeqn
Straightforward calculations give:
\[ \left(
\begin{array}{ccc}
b_{1;1} & b_{1;2} & b_{1;3} \\
b_{2;1} & b_{2;2} & b_{2;3} \\
b_{3;1} & b_{3;2} & b_{3;3}
\end{array}
\right)=\left(
\begin{array}{ccc}
0 & 0 & 0 \\
0 & 0 & -\delta \\
0 & \delta & 0
\end{array}
\right).
\]
Then, we have $r_{ij}=0$, that is, $\beta$ is a Killing 1-form with respect to $\alpha$. Further, we have
\[
\begin{array}{lll}
s_{11}=0, & s_{12}=0, & s_{13}=0, \\
s_{21}=0, & s_{22}=0, & s_{23}=-\delta, \\
s_{31}=0, & s_{32}=\delta, & s_{33}=0.
\end{array}
\]
It is obvious that $s_{1}=s_{2}=s_{3}=0$. Hence, we have
\be
s_{0}=0, \ \ s_{0;0}=0. \label{s00}
\ee
By the definition, $t_{ij}=s_{im}s^{m}_{ \ j}$, we can get
\[
\begin{array}{lll}
t_{11}=0, & t_{12}=0, & t_{13}=0, \\
t_{21}=0, & t_{22}=-\delta ^2, & t_{23}=0, \\
t_{31}=0, & t_{32}=0, & t_{33}=-\delta ^2.
\end{array}
\]
Then we have
\be
t_{00}=-\delta ^{2}\left((y^2)^{2}+(y^3)^{2}\right) \label{t00}
\ee
and $t_{1}=t_{2}=t_{3}=0$, which means
\be
t_{0}=0.  \label{t0}
\ee
Further, by  $s_{ij;k}=(ds_{ij}-s_{pj}\omega _{i}^{\ p}-s_{ip}\omega _{j}^{\ p})(e_{k})$, we have the following
\[
\begin{array}{lll}
s_{11;1}=0, & s_{12;1}=0, & s_{13;1}=0, \\
s_{21;2}=\delta\epsilon, & s_{22;2}=0, & s_{23;2}=0, \\
s_{31;3}=\delta\varepsilon, & s_{32;3}=0, & s_{33;3}=0.
\end{array}
\]
It is easy to see that
\be
s^{m}_{\ 0;m}= 2\delta\varepsilon y^1. \label{sm0m}
\ee

Then, from (\ref{alpha}), (\ref{beta}) and (\ref{s00}), (\ref{t00})  and (\ref{sm0m}), we find that (\ref{main11}) and (\ref{main12}) hold for
\be
c(x)=-(K-1). \label{c1}
\ee
At the same time, by (\ref{t0}), it is easy to check that (\ref{thm2}) holds for $\lambda = K$ and $c(x)=-(K-1)$.

\vskip 8mm

\noindent
Xinyue Cheng \\
School of Mathematical Sciences \\
Chongqing Normal University \\
Chongqing  401331,  P. R. of China  \\
E-mail: chengxy@cqut.edu.cn  \& chengxy@cqnu.edu.cn

\vskip 4mm

\noindent
Li Yin \\
School of Sciences \\
Chongqing University of Technology \\
Chongqing 400054,  P. R. China \\
E-mail: yinli@2015.cqut.edu.cn

\vskip 4mm

\noindent
Tingting Li \\
School of Sciences \\
Chongqing University of Technology \\
Chongqing 400054,  P. R. China \\
E-mail: litt@2015.cqut.edu.cn

\end{document}